\documentclass[12pt]{amsart}
\usepackage{amsmath,amssymb, euscript}

\theoremstyle{definition}

\begin{document}

\address{Azer Akhmedov, Department of Mathematics,
North Dakota State University,
Fargo, ND, 58102, USA}
\email{azer.akhmedov@ndsu.edu}
 
\begin{center} {\bf ON FREE DISCRETE SUBGROUPS OF \em{Diff}$(I)$} \end{center} 

\vspace{0.6cm}

\begin{center} {\Small AZER AKHMEDOV} \end{center}
 
 \vspace{1cm}
 
 {\Small ABSTRACT: We prove that a free group $\mathbb{F}_2$ admits a faithful discrete representation into {\em Diff}$_{+}^1[0,1]$. We also prove that $\mathbb{F}_2$ admits a faithful discrete representation of bi-Lipschitz class into {\em Homeo}$_{+}[0,1]$. Some properties of these representations have been studied. In the last section of the paper we raise several questions}. 
 
\vspace{1cm}
 
      \begin{center} INTRODUCTION \end{center}
      
      \medskip
      
  In recent decades and especially in recent years, some remarkable papers  have appeared which are devoted to the study of finitely generated subgroups of {\em Diff}$_{+}^1[0,1]$. In contrast, discrete subgroups of {\em Diff}$_{+}^1[0,1]$ is much less studied. Very little is known in this area especially in comparison with the very rich theory of discrete subgroups of Lie groups which has started in the works of F.Klein and H.Poincare in the 19th century, and has experienced enormous growth in the works of A.Selberg, A.Borel, G.Mostow, G.Margulis and many others in the 20th century. Many questions which are either very easy or have been studied long time ago for (discrete) subgroups of Lie groups remain open in the context of the infinite-dimensional group {\em Diff}$_{+}^1[0,1]$ and its relatives. In this paper, we address a question about existence of discrete faithful representation of non-abelian free groups into {\em Diff}$_{+}^1[0,1]$.
  
  \medskip
  
  We assume the usual topology on the group {\em Diff}$_{+}^1[0,1]$ given by the standard metric of $C^1[0,1]$. We will denote this metric by $d_1$.
   
  \medskip
  
  {\bf Theorem 1.} A free group $\mathbb{F}_2$ admits a faithful discrete representation into {\em Diff}$_{+}^1[0,1]$. 
  
  \bigskip
    
   We will also be interested in discrete subgroups of  {\em Homeo}$_{+}[0,1]$ - the group of orientation preserving homeomorphisms of the closed interval. Here, the metric comes from the sup norm of the Banach space $C[0,1]$. 
   
   \medskip
   
   {\bf Theorem 2.} A free group $\mathbb{F}_2$ admits a faithful discrete representation into {\em Homeo}$_{+}[0,1]$. Moreover, 
  
  \medskip
  
  a) the representation can be chosen from the class $C^1(0,1)\cap BiLip[0,1]$.
  
  \medskip
  
  b) for any non-empty open neighborhood $\Omega $ of {\em Homeo}$_{+}[0,1]$, the generators of the faithful discrete representation of $\mathbb{F}_2$ can be chosen from $\Omega $.
  
  \bigskip
  
  Here, $BiLip[0,1]$ denotes the set of all bi-Lipschitz functions from the closed interval $[0,1]$ into itself.  
     
  \vspace{0.7cm}
   
   \begin{center} PROOFS OF MAIN THEOREMS \end{center} 
   
   \medskip
   
   In this section we will prove Theorem 1 and Theorem 2.
   
   \medskip
   
   In the free group $\mathbb{F}_2$ we will fix the left-invariant Cayley metric with respect to standard generating set, and denote it by $|.|$. The following notions will be useful
   
   \medskip
   
   {\bf Definition 1.} Let $W$ be a reduced word in the alphabet of the standard generating set of the free group $\mathbb{F}_2$. We say that a reduced word $U$ is a {\em suffix} of $W$, if $W = U_1U$ where $U_1$ is a reduced word, and $|W| = |U_1| + |U|$. We also say that a reduced word $V$ is a {\em prefix} of $W$, if $W = VV_1$ where $V_1$ is a reduced word, and $|W| = |V| + |V_1|$. 
   
   \bigskip
   
   {\bf Proof of Theorem 1.} 
 
 \medskip
 
   Let $I_n = (\frac {1}{2n+1}, \frac {1}{2n})$ for any $n\in \mathbb{N}$ and let $C > 1$.
   
   \medskip
   
   We will build two maps $f,g\in ${\em Diff}$_{+}^1[0,1]$ such that the group $\Gamma _{f,g}$ generated by them is isomorphic to $\mathbb{F}_2$ and satisfies the following condition:
   
   $(c)$ for all $g_1, g_2\in \Gamma _{f,g}, g_1\neq g_2$, the inequality $sup_{t\in [0,1]} |g_1'(t) - g_2'(t)| > C$ is satisfied.
   
   \medskip
  
  Let $\pi _n = (U_n,V_n), n\geq 1$ be a sequence of pairs of words (elements) in $\mathbb{F}_2$ satisfying the following conditions:
  
   $(a1) |U_n| \geq |V_n|$ for all $n\geq 1$.
   
   \medskip
   
  $(a2)$ if $m > n$ then $|U_m| \geq |U_n|$.
    
   \medskip
   
   $(a3)$ if $m > n, |U_m| = |U_n|$ then $|V_m|\geq |V_n|$.
   
   \medskip
   
   $(a4)$ $U_n\neq 1$, for all $n\geq 1$.
   
   \medskip
   
   $(a5)$ if $U, V\in \mathbb{F}_2, U\neq 1, |U|\geq |V|$ then there exists $n\in \mathbb{N}$ such that $U = U_n, V = V_n$.
   
   \medskip
   
   $(a6)$ if $m\neq n$ then $\pi _m\neq \pi_n$.
   
   \medskip
   
   For every $n$, the longest common suffix of $U_n$ and $V_n$ will be denoted by $W_n$ and $s_n = |W_n|$.
   
   \medskip
   
  Let also $m_n = Card \{n \ | \ \pi_n = (U_n, V_n), |U_n| = n\}$ for all $n\geq 1, m_0 = 0$. So, indeed $m_n, n\geq 1$ is the number of words of length $n$ (which equals to $4\times 3^{n-1}$).
  
  \medskip
  
  Let $\alpha = (\alpha _1, \alpha _2, \ldots )$ be a sequence of positive real numbers such that
  
  \medskip
  
  $(b1)$ $\lim _{r\rightarrow \infty } \alpha _r = 0$
  
  \medskip
  
  $(b2)$ for every $r\in \mathbb{N}, s\in \{0, 1, \ldots , r-1\}$, the inequality $$(1+\alpha _r)^s[(1+\alpha _r)^{r-s}-1] > C$$ is satisfied.
       
  \medskip
  
  Let also $\beta = (\beta _1, \beta _2, \ldots )$ be a sequence such that $\beta _i = \alpha _j$ for all $m_1 + \ldots m_{j-1} < i \leq m_1 + \ldots m_{j-1} + m_j$. We notice that $\lim _{n\rightarrow \infty }\beta _n = 0$; moreover, for every $n\in \mathbb{N}$, $\beta _n = \alpha _{i(n)}$ where $i(n)\rightarrow \infty $ as $n\rightarrow \infty $.
  
  \medskip
  
 Now, for any natural $n$, let $x_0^n$ be the midpoint of the interval $I_n$, $s = s_n$, and let $f, g$ be defined in the interval $I_n$ such that 
 
 \medskip
 
  $(c1) f(x) = g(x) = x$ for all $x\in \{\frac {1}{2n+1}, \frac {1}{2n}\}$
  
  \medskip
  
  $(c2) f'(x)\in [\frac {1}{1+ \beta _n + \frac{1}{n}}, 1 + \beta _n + \frac{1}{n}]$, for all $x\in I_n$.
   
   \medskip
   
  $(c3) f'(x) = g'(x) = 1$ for all $x\in \{\frac {1}{2n+1}, \frac {1}{2n}\}$
   
   \medskip
   
  $(c4)$ if $|U_n| = r, U_n = a_ra_{r-1}\ldots a_s\ldots a_1, a_i\in \{f, g, f^{-1}, g^{-1}\}, 1\leq i\leq n, U_n(k) = a_k\ldots a_1, \ 0\leq k\leq n$ then $a_{k+1}'(U_n(k)(x_0^n)) = 1 + \beta _n$. 
  
  \medskip
  
  $(c5)$ if $|V_n| = m, \ V_n = b_mb_{m-1}\ldots b_1, b_i\in \{f, g, f^{-1}, g^{-1}\}, 1\leq i\leq m$, \ and if $m-1 \geq k\geq s, V_n(k) = b_k\ldots b_1$ then the equality $b_{k+1}'(V_n(k)(x_0^n)) = 1$ holds.
  
  \medskip
  
  Now, if $x\in [0,1]\backslash \sqcup _{n\in \mathbb{N}}I_n$, we set $f(x) = g(x) = x$ (hence  $f'(x) = g'(x) = 1$).
  
  \medskip

 Then the functions $f,g$ will belong to {\em Diff}$^1[0,1]$. Moreover, for any $n\geq 1$, by Chain Rule, we have  $$U_n'(x_0^n) = (1+\beta _n)^r, V_n'(x_0^n) = (1+\beta _n)^{s} 1^{m-s} = (1+\beta _n)^{s}$$    
 
 \medskip

  Since $\beta _n = \alpha _{i(n)}$ and $i(n) = r$, the inequality $|(U_n(f,g))'(x_0^n) - (V_n(f,g))'(x_0^n)| > C$ follows from condition $(b2)$. $\square $
  
 \bigskip
 
 {\bf Remark 1.} We indeed prove more than discreteness; the inequality $$\sup _{t\in [0,1], g\in \mathbb{F}_2\backslash \{1\}}|g'(t)-1| \geq  C > 0$$ would suffice for discreteness. By proving more general inequality $$\sup _{t\in [0,1], g_1, g_2\in \mathbb{F}_2, g_1\neq g_2}|g_1'(t)-g_2'(t)| \geq  C > 0$$ we show that the representation is {\em uniformly discrete.} Since the metric in {\em Diff}$_{+}^1[0,1]$ is not left-invariant, discreteness does not necessarily imply uniform discreteness.  
 
 \bigskip
 
 {\bf Remark 2.} It is clear from the proof that the functions $f(t)$ and $g(t)$ can be chosen from arbitrary non-empty open neighborhood of identity. This is contrary to the case of connected Lie groups: {\em Margulis Lemma} states that any connected Lie group $G$ possesses a non-empty open neighborhood $U$ of identity such that any discrete subgroup of $G$ generated by elements from $U$ is nilpotent (see {\bf [R]}). Thus we have shown that Margulis Lemma does not hold for the group {\em Diff}$^1[0,1]$.
 
 \bigskip
 
  It is easy to put the proof (main idea of the proof) of Theorem 1 in words: we take all pairs $(U_n, V_n)$ in the free group $\mathbb{F}_2$ that is interesting to us and enumarte them with some care, [conditions $(a1)-(a6)$]. Then we choose countable pairwise disjoint open subintervals $I_1, I_2, \ldots I_n, \ldots $ of [0,1] which are accumulating to the left endpoint of $[0,1]$, ($I_i$ is on the left side of $I_j$ for all $i > j$). Then at each of subintervals we arrange the maps $f,g$ such that  $sup_{x\in I_n} |(U_n(f,g))'(x)-1|$ and $sup_{x\in I_n} |(V_n(f,g))'(x)-1|$ converge to zero as $n\rightarrow \infty $, moreover, at some point $x_n\in I_n$ we have $$U_n'(x_n) > C, \ V_n'(x_n) = 1$$ Because of decay of derivatives to 1, the second condtion (the equality) becomes a natural condition. As for the first condition (the inequality), one notices that the word $U_n$ has length at least $log(n)$ which goes to infinity as $n$ grows. Then, since  $U_n'(x_n)$ is the product of $log(n)$ derivatives we can have this product to be bigger than $C$ yet each of the factor stay close to 1. (and converge to 1 as $n$ goes to infinity). For fixed $n$, each of these conditions imposes only finitely many conditions on $f$ and $g$ in $I_n$, and for the next pair we go to a different interval $I_{n+1}$, hence we have no obstruction left to the existence of discrete $F_2$ of $C^1$ class.
     
 \medskip
 
 However, because of slow growth of $log(n)$, and because the lengths of intervals of $I_n$ converge to zero faster than $1/n$, it is easy to see that this construction will not work in $C^2$ class, in fact, as Danny Calegari pointed out, it will not work in any $C^{1+\epsilon }$ class for any $\epsilon > 0$; imposing same condition will blow-up the Holder norm. So one cannot achieve higher regularity of representations by taking care of different pairs in disjoint areas of the closed interval $[0,1]$. If we want to mix fields of actions for different pairs, we need to take some cautions.        
 
 \bigskip

 Now we will prove Theorem 2. We need the following definitions
 
 \medskip
 
 {\bf Definition 2}. For open subintervals $I,J\subset (0,1)$ we say $I<J$ if any element  $I$ is less than any element of $J$.
 
 \medskip
  
 {\bf Definition 3.} A two-sided sequence $\{I_n\}_{n\in \mathbb{Z}}$ of open subintervals of $(0,1)$ is called {\em a chain} if $I_n < I_{n+1}$ for all $n\in \mathbb{Z}$. 

  \bigskip

 {\bf Proof of Theorem 2.} Let $\epsilon > 0$, and let $A_n, B_n, \ n\in \mathbb{Z}$ be open subintervals of $(0,1)$ such that 
 
  \medskip

 (i) the two sided sequence $\{A_n, B_n\}_{n\in \mathbb{Z}}$ is a chain of subintervals.

  [i.e. we have $\ldots < A_{-1} < B_{-1} < A_{0} < B_{0} < A_{1} < B_{1} < A_{2} < \ldots $]
 
 \medskip
 
 (ii) for all $n\in \mathbb{Z}$, and for all $i\in \{1,2,3,4\}$ we have
 
     $f^i(A_{n})\subseteq B_{n}, f^{-i}(A_{n})\subseteq B_{n-1}$.
     
 \medskip
     
  (iii)   for all $n\in \mathbb{Z}$,
   
     $$g(B_{n})\subseteq \cup _{n\in Z}A_n, g^{-1}(B_{n})\subseteq \cup _{n\in Z}A_n$$
  
 \medskip
 
 (iv) for all $n\in \mathbb{Z}$, the inequality \ $sup_{x\in A_n, y\in A_{n+2}}|x-y| < \epsilon $ holds. 
    
   \bigskip
      
 It is starightforward to choose $f,g\in ${\em Homeo}$_{+}[0,1]$ satisfying conditions (i)-(iv).

\medskip

 Now, let $A = \cup _{n\in Z}A_n, B = \cup _{n\in Z}B_n$. Notice that by conditions (i)-(ii)

 $$f^i(A)\subseteq B  \ \forall i\in \{-4, -3, -2, -1, 1, 2, 3, 4\}, \  g^i(B)\subseteq A, \ \forall i\in \{-1, 1\}$$
  
  This allows us to use ping-pong argument.

\medskip

 Ping-pong argument is usually used to guarantee existence of free subgroups, here we will be using it also to satisfy discreteness (which is natural). Using ping-pong lemma, we will show that, assuming conditions (i)-(iv), if $U(f,g), V(f,g)$ are reduced words satisfying conditions

  (1) $U(f,g) = f^2U_0(f,g)f^2, \ V(f,g) = fV_0(f,g)f$ where $U_0(f,g), V_0(f,g)$ are both non-empty reduced words starting and ending in letter $g$
  
  (2) $|U_0(f,g)| = |U(f,g)| - 4,  |V_0(f,g)| = |V(f,g)| - 2$.
    
 (3) none of the letters $\{f, g, f^{-1}, g^{-1}\}$ occur with exponent other than $\{-1,1\}$ in $U_0(f,g)$ and in $V_0(f,g)$.
 
 \medskip
 
 then $U(f,g), V(f,g)$ actually generate a free subgroup isomorphic to $\mathbb{F}_2$ in {\em Homeo}$_{+}[0,1]$. We will have that this subgroup (which we will denote by $\Gamma $) is discrete.   

\medskip

  Let $W(U,V)$ be any reduced non-trivial word in the alphabet \ $\{U = U(f,g), V = V(f,g), U^{-1} = U(f,g)^{-1}, V^{-1} = V(f,g)^{-1}\}$. Then in the alphabet $\{f,g, f^{-1}, g^{-1}\}$ the word $W$ ends with either $f$ or $f^{-1}$. 
 
 \medskip
 
 Let  $x_0$ be the midpoint of $A_0$. 
 
 \medskip
 
 We notice that $f^i(A)\subseteq B$ for all $i\in \{-4, -3, -2, -1, 1, 2, 3, 4\}$. Furthermore, $g^{\pm 1}(B)\subseteq A$. Then by standard ping-pong argument, we have that $$W(x_0) = W(U(f,g),V(f,g))(x_0)\notin A_0$$ hence $W\neq 1$ in $\Gamma $, and $|W|_0 \geq |A_0|/2$. 
 
 \medskip
 
 We now consider the general case of arbitrary distinct $h_1, h_2\in \Gamma $.
Let $W_1(U,V), W_2(U,V)$ be two distinct reduced words in the alphabet $\{U, V, U^{-1}, V^{-1}\}$.  Then we can write $W_2 = WW_1$ where $W = W(U,V) = W(U(f,g), V(f,g))$ satisfies conditions similar to (1)-(3).

\medskip

 Since $W_1 = W_1(U(f,g), V(f,g))$ is a bijective map from $[0,1]$ onto
$[0,1]$, there exists $z\in [0,1]$ such that $W_1(z) = x_0$. Then $W_2(z) = W(W_1(z)) = W(x_0) \notin A_0$. 

 \medskip

 Thus we have $|W_1(z) - W_2(z)| = |x_0 - W(x_0)| > |A_0|/2$. Thus we established that the non-abelain free subgroup generated by $U$ and $V$ is discrete.
 
 \medskip
 
 For the claim b), let $\Omega $ contains a ball of radius $r$, and $M = max \{|U|,|V|\}$. Then by condition (iv), $max\{||U||_0, ||V||_0\} < \epsilon M$. We can choose $U, V$ such that $M < 100$. Since $\epsilon $ is arbitrary we can choose it to be such that $100\epsilon  < r$, and hence we obtain the claim b).
 
 \medskip
 
 For the claim a), we may further assume that $$A_n = (\frac{1}{5(|n|+1)}, \frac{1}{5|n|+4}), B_n =  (\frac{1}{5|n|+4}, \frac{1}{5|n|)}), \ \forall n\leq -2 $$ and $$A_n = (1 - \frac{1}{5n}, 1 - \frac{1}{5n+1}), B_n =  (1 - \frac{1}{5n+1}, 1 - \frac{1}{5(n+1)}), \ \forall n\geq 1$$ (and we choose $A_{-1}, B_{-1}, A_0, B_0$ to be arbitrary open non-empty intervals such that condition (i) holds). Then it is straightforward to choose $f, g\in ${\em Homeo}$_{+}[0,1]$ such that $f\in C^{1}[0,1], g\in C^1(0,1),$ and $g$ is a  bi-Lipschitz function with Lipschitz constant at most 5 in $[0,\frac{1}{5}]$ and in $[\frac{4}{5},1]$. Then for any word $W$ in the free group $\mathbb{F}_2$, the function $W(U(f,g), V(f,g))$ will be a bi-Lipshitz function of class $C^1(0,1)$.$\square $

 \bigskip
    
 {\bf Remark 3.} We would like to point out what goes wrong if one applies the idea of the proof to Theorem 2 directly to obtain a faithful discrete representation of $\mathbb{F}_2$ in {\em Diff}$_{+}^1[0,1]$ :  
  
  \medskip
  
 Let $A_n,B_n, n\in \mathbb{Z}$ be mutually disjoint open subintervals in $(0, 1)$ satisfying conditions (i), (ii) and (iii).
 
 \medskip
 
  We will show that it is impossible to have the maps differentiable ($C^1$ class) under these conditions (i)-(iii); there are obstructions easily obtained from Mean Value Theorem. Indeed, let $\lim _{x\rightarrow 1-}f(x) = p$. (Then $p\leq 1$). 
  
  \medskip
  
  Let $p_1, p_2$ be positive real numbers such that $$p_1 < p < p_2, \ p_1 > \frac{99}{100}p, \ p_2 < \frac{101}{100}p$$
  
  \medskip

  Then, by Mean Value Theorem, from condition (ii)$'$ we obtain that $$|B_n| >(p_1 + p_1^2 + p_1^3)|A_n| \ and  \ |B_n| > (\frac {1}{p_2} + \frac {1}{p_2^2} + \frac {1}{p_2^3})|A_{n+1}|$$ for sufficiently big positive $n$. Then $$\frac {|g(B_n)|}{|B_n|} \leq \frac {|A_{n+1}|}{|B_n|} < \frac {1}{\frac{1}{p_2} + \frac{1}{p_2^2} + \frac{1}{p_2^3}}$$ and $$\frac {|g^{-1}(B_n)|}{|B_n|} \leq \frac {|A_n|}{|B_n|} < \frac{1}{p_1 + p_1^2 + p_1^3}$$
  
  \medskip
  
  Then, by Mean Value Theorem, we obtain that for sufficiently big positive $n$, there exists $u_n, v_n\in B_n$ such that $g'(u_n) < \frac {1}{\frac{1}{p_2} + \frac{1}{p_2^2} + \frac{1}{p_2^3}}$, and $(g^{-1})'(v_n) < \frac{1}{p_1 + p_1^2 + p_1^3}$.  
  
   \medskip
   
  However, since $\lim_{x\rightarrow 1-} g'(x) = 1/\lim_{x\rightarrow 1-} (g^{-1})'(x)$, we obtain a contradiction because   $$\frac {1}{\frac{1}{p_2} + \frac{1}{p_2^2} + \frac{1}{p_2^3}} \ \frac{1}{p_1 + p_1^2 + p_1^3} < \frac{1}{\frac{p_1}{p_2} + \frac{p_1^2}{p_2^2} + \frac{p_1^3}{p_2^3}} < \frac {1}{2}   < 1$$

  \medskip
  
  {\bf Remark 4.} In the proof of Theorem 2, by slightly changing conditions (1)-(3), it is possible to replace condition (ii) by the following weaker version:
    
  (ii)$'$ for all $i\in \{1,2\}, n\in \mathbb{Z}$, \ $f^i(A_n)\subseteq B_n, f^{-i}(A_n)\subseteq B_{n-1}$
  
  \medskip
  
  However, a similar argument shows that there are no $f, g\in ${\em Diff}$_{+}^1[0,1]$ satisfying conditions (i), (ii)$'$ and (iii). It also follows from the criterion of D.Calegari in {\bf [C]} that no $C^1$-class diffeomorphisms exist which satisfy conditions (i), (ii)$'$ and (iii).
  
 \bigskip
 
 {\bf Remark 5.} The metric in $C^1[0,1]$ is given by the norm $||f|| = ||f||_0 + ||f||_1$ where $||f||_0 = sup _{x\in [0,1]}|f(x)|, ||f||_1 = sup _{x\in [0,1]}|f'(x)|$. If $||f||_1$ is small and $f(0) = 0$, then by Mean Value Theorem $||f||_0$ cannot be big. However, $||f||_1$ can be big even if $||f||_0$ is small. In the proof of Theorem 1, taking $f(x) = W(x) - x$, we actually show that $||f||_1$ stays big for all $W\neq 1$. We do not show that $||f||_0$ is big, in fact, $||f||_0$ converges to zero as we go through the sequence $W$ in the order as described in the proof. However, in the proof of Theorem 2, we indeed show a stronger fact that $||f||_0$ remains big.    

 \bigskip
   
  \vspace{0.7cm}
 
 \begin{center} QUESTIONS \end{center}
 
  \bigskip
  
  In this section, we raise several questions. We will address these questions in our next publication.
  
  \medskip
  
  The regularity of the represenation is a very interesting question; if a finitely generated group $\Gamma $ admits a faithful discrete representation in {\em Diff}$_{+}[0,1]$ or in {\em Homeo}$_{+}[0,1]$, it is interesting to know if one can achieve faithful discrete representations of higher ($C^k, k > 1, \ C^{\infty }$, analytic, etc)  regularity. 
  
   \medskip
   
   {\bf Question 1.} Does a free group $\mathbb{F}_2$ admit a faithful discrete representation into {\em Diff}$_{+}[0,1]$ 
   
   a) of $C^k$ regularity for some $k>1$? 
  
   b) of $C^k$ regularity for any $k\geq1$? 
  
   c) of $C^{\infty }$ regularity?
   
   \bigskip
   
   {\bf Question 2.} Does $\mathbb{F}_2$ admit a faithful discrete analytic representation into {\em Diff}$_{+}[0,1]$? 
 
   \bigskip
   
   Let $\Gamma $ be a finitely generated group, and $\pi :\Gamma \rightarrow ${\em Diff}$_{+}^1[0,1]$ be a faithful discrete representation of it.

   \medskip
   
   {\bf Definition 4.} The representation $\pi $ is called {\em $||.||_0$-discrete} if there exists $C>0$ such that $||\pi (g)||_0 > C$ for all $g\in \Gamma \backslash \{1\}$. 
 
 \medskip
 
  By Remark 5, $||.||_0$-discreteness of the representation implies its discreteness. Also, $||.||_0$-discrete representation of a group into {\em Diff}$_{+}^1[0,1]$ is just a discrete representation into {\em Homeo}$_{+}[0,1]$ of $C^1$-regularity.   

 \medskip
 
  {\bf Question 3.} Does $\mathbb{F}_2$ admit a faithful $||.||_0$-discrete representation into {\em Diff}$_{+}^1[0,1]$? 
  
 \bigskip
  
   {\bf Definition 5.} The representation $\pi $ is called {\em strongly discrete} if there exists $C>0$ and $x_0\in (0,1)$ such that $||\pi(g)(x_0)||_1 > C$ for all $g\in \Gamma \backslash \{1\}$. 
  
  \medskip
  
  {\bf Question 4.} Does $\mathbb{F}_2$ admit a faithful strongly discrete representation into {\em Diff}$_{+}^1[0,1]$?
  
  \bigskip
  
  Similarly, we say that a faithful representation $\pi :\Gamma \rightarrow ${\em Homeo}$_{+}[0,1]$ is {\em strongly discrete} (in {\em Homeo}$_{+}[0,1]$) if there exists $C>0$ and $x_0\in (0,1)$ such that $||\pi(g)(x_0)||_0 > C$ for all $g\in \Gamma \backslash \{1\}$. Notice that in the proof of Theorem 2, the representation of $\mathbb{F}_2$ into {\em Homeo}$_{+}[0,1]$ is indeed strongly discrete.
    
  \bigskip

  {\bf Definition 6.} Let $G$ be a topological group or a group with a metric. We say that {\em Weak Margulis Lemma } holds for $G$, if there exists an open non-empty neighborhood $U$ of identity such that any discrete subgroup of $G$ generated by elements from $U$ does not contain a non-abelian free subgroup. 
  
  \medskip
   
   We will be interested in the group {\em Diff}$_{+}^{1+\epsilon }[0,1]$ where $\epsilon $ is a fixed positive real number. On this group, we are considering the metric $d_1$, i.e. the metric which comes from the Banach norm of $C^1[0,1]$
   
   \medskip
    
   {\bf Question 5.} Does Weak Margulis Lemma hold for the group {\em Diff}$_{+}^{1+\epsilon }[0,1]$ for some $\epsilon >0$?  
   
   \medskip
    
   {\bf Remark 6.} It follows from the proof of Theorem 1 and from Theorem 2 that Weak Margulis Lemma does not hold neither for {\em Diff}$_{+}^{1}[0,1]$ nor for {\em Homeo}$_{+}[0,1]$, in respective metrics.
   
   \medskip
   
  The study of discrete subgroups of {\em Diff}$_{+}^{1}[0,1]$ is interesting beyond the existence question of discrete faithful represenations of free groups or even of the groups which contain non-abelian free subgroups. The existence of discrete faithful representation into {\em Diff}$_{+}^{1}[0,1]$ imposes some algebraic properties onto the group; for example, it is well-known that any subgroup of {\em Homeo}$_{+}[0,1]$ is left-orderable. Furthermore, if a group is isomorphic to a subgroup of {\em Diff}$_{+}^{1}[0,1]$ then it is locally indicable, as proven by W.Thurston ({\bf [T]}). It is interesting if discreteness implies further algebraic restrictions on the group. We would like to ask the following  
  
  \medskip
  
   {\bf Question 7.} Is there a finitely generated group which admits a faithful representation into {\em Diff}$_{+}^{1}[0,1]$ but does not admit a faithful discrete representation?
   
   \bigskip
   
   {\em Acknowledgment:} I am thankful to Matthew G.Brin and Danny Calegari for useful discussions related to the content of this paper.
      
     \vspace{2cm}
     
     {\bf R e f e r e n c e s:}
     
     \bigskip
     
     [C] Calegari, D.  \ Nonsmoothable, locally indicable group  
actions on the interval. \ {\em Algebraic and Geometric Topology} {\bf 8} (2008) no. 1, 609-613     

  \medskip

     [R] Raghunathan,M.S.  \ Discrete subgroups of semi-simple Lie groups. \ Band 68, Springer-Verlag, New York-Heidelberg, 1972.

   \medskip
        
     [T]  Thurston, W. \ A generalization of Reeb stability theorem. \ {\em Topology} {\bf 13}, (1974) 347-352. 
     
\end{document}